\documentclass[12pt]{amsart}
\usepackage{amsmath,amssymb,bm,mathrsfs}
\usepackage{graphicx,color}
\usepackage{subfigure}
\usepackage{tikz}
\newtheorem{theorem}{Theorem}[section]
\newtheorem{lemma}[theorem]{Lemma}
\newtheorem{remark}[theorem]{Remark}
\numberwithin{equation}{section}
\allowdisplaybreaks[4]
\begin{document}
\def\O{\Omega}
\def\R{\mathbb{R}}
\def\p{\partial}
\def\cT{\mathcal{T}}
\def\LT{{L_2(\O)}}
\def\HOne{{H^1(\O)}}
\def\HOnez{H^1_0(\O)}
\def\Uad{U_{\rm ad}}
\def\argmin{\mathop{\rm argmin}}
\def\CPF{\mathrm{C}_{\rm PF}}
\def\bysd{\bar y_{*,\dag}}
\def\bpsd{\bar p_{*,\dag}}
\def\busd{\bar u_{*,\dag}}
\def\tud{\tilde u_\dag}
\def\tys{\tilde y_*}
\def\Ksd{K_{*,\dag}}
\def\DES{V_*\times W_\dag}
\def\sss{\scriptscriptstyle}
\def\red{\color{red}}
\title[New Estimates for An Elliptic Optimal Control Problem]
{New Error Estimates for An Elliptic Distributed Optimal Control
Problem with Pointwise Control Constraints}
\author{Susanne C. Brenner}
\address{S.C. Brenner, Department of Mathematics and Center for Computation \& Technology,
 Louisiana State University, Baton Rouge, LA 70803, USA}
\email{brenner@math.lsu.edu}
\author{Li-yeng Sung}
\address{L.-Y. Sung, Department of Mathematics and Center for Computation \& Technology,
 Louisiana State University, Baton Rouge, LA 70803, USA}
\email{sung@math.lsu.edu}
\thanks{This work  was supported in part
 by the National Science Foundation under
 Grant No. DMS-22-08404, and by the National Science Foundation under Grant No. DMS-1929284 while
 the authors were in residence at  the Institute for Computational and Experimental Research in Mathematics in Providence, RI, during the semester program in Spring 2024.}
\begin{abstract} We derive error estimates for a linear-quadratic
elliptic distributed
optimal control problem with pointwise control constraints that
can be applied to
standard finite element methods and multiscale finite element methods.
\end{abstract}
\keywords{elliptic distributed optimal control problems, pointwise
control constraints,
 finite element methods, rough coefficients, multiscale}
\subjclass{65N30, 65K10, 49M41}
\date{\today}
\maketitle
\section{Introduction}\label{sec:Introduction}
 Let $\O$ be a polygonal (resp., polyhedral) domain in $\R^d$ for $d=2$
 (resp., $3$), $y_d\in\LT$,
 $f\in\LT$ and $\gamma\leq 1$ be
 a positive constant.   The optimal control problem is to find
\begin{equation}\label{eq:OCP}
   (\bar y,\bar u)= \argmin_{(y,u)\in K}\frac12\Big[\|y-y_d\|_\LT^2
   +\gamma\|u\|_\LT^2\Big],
\end{equation}
 where $(y,u)$ belongs to $K\subset \HOnez\times\LT$ if and only if
\begin{equation}\label{eq:PDEConstraint}
  a(y,z)=\int_\O (f+u)z\,dx\qquad \forall\, z\in\HOnez
\end{equation}
 and
\begin{equation}\label{eq:ControlConstraints}
  u\in \Uad=\{v\in\LT:\,\phi_1\leq v\leq \phi_2 \;\;\text{in}\;\;\O\}.
\end{equation}
\par
 Here the symmetric bilinear form $a(\cdot,\cdot)$ on $\HOne$ satisfies
\begin{equation}\label{eq:BilinearForm}
  \alpha|y|_\HOne^2\leq a(y,y) \leq \beta |y|_\HOne^2
  \qquad\forall\,y\in\HOne,
\end{equation}
 where $\alpha\leq\beta$ are positive constants,
 and we assume that
\begin{equation}\label{eq:ConstraintsRegulaity}
  \phi_1,\phi_2\in\HOne
\end{equation}
 satisfy
\begin{equation}\label{eq:ControlCompatibility}
  \phi_1\leq\phi_2 \;\;\text{in}\;\;\O.
\end{equation}
\begin{remark}\label{rem:AlmostEverywhere}\rm
 Throughout this paper
  the inequalities and equalities between functions  are to be interpreted
   in the sense of almost everywhere in $\O$.
\end{remark}
\begin{remark}\label{rem:Notation}\rm
  We follow the standard notation for function spaces, norms and differential
   operators that can be found for example
  in \cite{ADAMS:2003:Sobolev,BScott:2008:FEM}.
\end{remark}
\begin{remark}\label{rem:RC}\rm
  The condition \eqref{eq:BilinearForm} is satisfied by
  many partial differential equation
  constraints with rough coefficients.
\end{remark}
\par
 The optimal control problem defined by \eqref{eq:OCP}--\eqref{eq:BilinearForm}
  is a model linear-quadratic
 problem (cf. \cite{Lions:1971:OC,Troltzsch:2010:OC}) and the error analysis
  of a finite element method for this problem
 was first given in \cite{Falk:1973:Control} 
 under additional assumptions on the
 bilinear form $a(\cdot,\cdot)$.  A substantial literature has been developed over the years
 (cf. the monographs
 \cite{NST:2006:EllipticOptimalControl,LY:2008:Adaptive,HPUU:2009:Book} and the references therein).
 Nevertheless, the existing error analysis cannot be directly applied to multiscale finite element methods
 under the rough coefficient assumption in \eqref{eq:BilinearForm}.
 \par
 Our goal is to develop new abstract error estimates under the assumption
 \eqref{eq:BilinearForm} that are
 suitable for the error analysis of classical finite element methods and also for
 multiscale finite element methods.  Our results
  (cf. Theorem~\ref{thm:AbstractLTError} and Theorem~\ref{thm:AbstractEnergyError})
  reduce the error analysis of finite element methods for the optimal control problem to the error analysis of
  finite element methods for elliptic boundary value problems.  Therefore
  they can be applied to any finite element methods
  that have already been analysed for elliptic boundary value problems.  In particular they can be
  applied to many multiscale finite element methods.
\par
 The rest of the paper is organized as follows. We recall the relevant
 properties of the optimal
 control problem in Section~\ref{sec:Continuous} and introduce the approximation
 problem in Section~\ref{sec:Approximation}.  We derive the abstract error estimates
 in Section~\ref{sec:Analysis}, present several applications in
 Section~\ref{sec:Applications}
 and end with some concluding remarks in Section~\ref{sec:Conclusions}.
\section{The Continuous Problem}\label{sec:Continuous}
 According to the classical theory in
 \cite{ET:1999:Convex,KS:1980:VarInequalities},
 the convex minimization problem defined by
  \eqref{eq:OCP}--\eqref{eq:BilinearForm} and
 \eqref{eq:ControlCompatibility} has a unique solution
 $(\bar y,\bar u)\in K$
 characterized by the first order optimality condition
\begin{equation}\label{eq:OptimalityCondition}
  \int_\O (\bar y-y_d)(y-\bar y)dx+\gamma\int_\O \bar u(u-\bar u)dx
  \geq0\qquad\forall\,(y,u)\in K.
\end{equation}
\par
 Let the adjoint state $\bar p\in \HOnez$ be defined by
\begin{equation}\label{eq:AdjointState}
   a(q,\bar p)=\int_\O (\bar y-y_d) q\,dx \qquad\forall\,q\in\HOnez.
\end{equation}
 In view of \eqref{eq:PDEConstraint} and \eqref{eq:AdjointState},
 we have, for any $(y,u)\in K$,
\begin{align}\label{eq:SAS}
  \int_\O (\bar y-y_d)(y-\bar y)dx+\gamma\int_\O \bar u(u-\bar u)dx
  &=a(\bar y-y,\bar p)+\gamma\int_\O \bar u(u-\bar u)dx\\
  &=\int_\O \bar p(u-\bar u)dx+\gamma\int_\O \bar u(u-\bar u)dx,\notag
\end{align}
 and hence
\begin{equation*}
  \int_\O (\bar p+\gamma\bar u)(u-\bar u)dx\geq0 \qquad\forall\,u\in \Uad,
\end{equation*}
 which means that $\bar u$ is the $L_2$ projection of the function
 $-(\bar p/\gamma)$ on the closed convex subset $\Uad$ of
 $\LT$.
 Consequently we have
\begin{equation}\label{eq:baruAndbarp}
  \bar u=\max(\phi_1,\min(\phi_2,-(\bar p/\gamma)))
\end{equation}
 and $(\bar y,\bar u)$ is determined by \eqref{eq:PDEConstraint},
 \eqref{eq:AdjointState} and \eqref{eq:baruAndbarp}.
%
\subsection{Bounds for $\|\bar y-y_d\|_\LT$ and $\|\bar u\|_\LT$}
\label{subsec:OptimalStateand ControlBdd}
 It follows from \eqref{eq:OCP}, \eqref{eq:ControlConstraints},
 \eqref{eq:ControlCompatibility} and $\gamma\leq1$ that
\begin{equation}\label{eq:AprioribaryEstimate}
  \|\bar y-y_d\|_\LT^2+\gamma\|\bar u\|_\LT^2
  \leq \|y_1-y_d\|_\LT^2+\|\phi_1\|_\LT^2,
\end{equation}
 where $y_1\in\HOnez$ is defined by
\begin{equation}\label{eq:y1}
  a(y_1,z)=\int_\O (f+\phi_1)z\,dx\qquad\forall\,z\in\HOnez.
\end{equation}
\par
From \eqref{eq:BilinearForm} and \eqref{eq:y1} we have
\begin{align*}
  \alpha|y_1|_\HOne^2\leq a(y_1,y_1)
        =\int_\O (f+\phi_1)y_1\,dx
        \leq \|f+\phi_1\|_\LT\|y_1\|_\LT,
\end{align*}
 which together with the Poincar\'e-Friedrichs inequality
\begin{equation}\label{eq:PF}
  \|v\|_\LT\leq \CPF|v|_\HOne \quad\forall\,v\in\HOnez
\end{equation}
 implies
\begin{equation}\label{eq:y1Bdd}
  \|y_1\|_\LT\leq (\CPF/\alpha)\|f+\phi_1\|_\LT.
\end{equation}
 Combining \eqref{eq:AprioribaryEstimate}, \eqref{eq:y1Bdd}
 and the Cauchy-Schwarz inequality,
 we find
\begin{equation*}
  \|\bar y-y_d\|_\LT^2+\gamma\|\bar u\|_\LT^2\leq 2\|y_d\|_\LT^2
  +4(\CPF^2/\alpha^2)\|f\|_\LT^2%
  +[4(\CPF^2/\alpha^2)+1]\|\phi_1\|_\LT^2.
\end{equation*}
\par
 Similarly we have
\begin{equation*}
  \|\bar y-y_d\|_\LT^2+\gamma\|\bar u\|_\LT^2\leq 2\|y_d\|_\LT^2
  +4(\CPF^2/\alpha^2)\|f\|_\LT^2%
  +[4(\CPF^2/\alpha^2)+1]\|\phi_2\|_\LT^2
\end{equation*}
 and hence
\begin{align}
   \|\bar y-y_d\|_\LT&\leq C_\sharp,\label{eq:OptimalStateBdd}\\
   \|\bar u\|_\LT&\leq \gamma^{-1} C_\sharp, \label{eq:baruLTBdd}
\end{align}
 where
\begin{equation}\label{eq:CSharp}
  C_\sharp=\Big(2\|y_d\|_\LT^2+4(\CPF^2/\alpha^2)\|f\|_\LT^2
  +[4(\CPF^2/\alpha^2)+1]\min(\|\phi_1\|_\LT^2,\|\phi_2\|_\LT^2)\Big)^\frac12.
\end{equation}
%
\subsection{Bounds for $|\bar u|_\HOne$ and $|\bar p|_\HOne$}\label{subsec:Bdds}
 It follows from \eqref{eq:BilinearForm} and \eqref{eq:AdjointState} that
\begin{equation*}
  \alpha|\bar p|_\HOne^2\leq a(\bar p,\bar p)=\int_\O (\bar y-y_d)\bar p\,dx\leq
   \|\bar y-y_d\|_\LT\|\bar p\|_\LT,
\end{equation*}
 which together with \eqref{eq:PF} and \eqref{eq:OptimalStateBdd} implies
\begin{equation}\label{eq:AdjointStateBdd}
  |\bar p|_\HOne\leq (\CPF/\alpha) C_\sharp.
\end{equation}
\par
 Since the space $\HOne$ is invariant under the max and min operators
 (cf. \cite[Lemma~7.6]{GT:2001:EllipticPDE}),
 we conclude from \eqref{eq:ConstraintsRegulaity} and \eqref{eq:baruAndbarp}
  that $\bar u\in \HOne$ and
\begin{equation}\label{eq:OptimalControlBdd}
  |\bar u|_\HOne\leq\max\big(|\phi_1|_\HOne,|\phi_2|_\HOne,\gamma^{-1}|\bar p|_\HOne\big).
\end{equation}
%
\subsection{The Lagrange Multiplier $\lambda$}\label{subsec:lambda}
 The function
\begin{equation}\label{eq:Multiplier}
  \lambda=\bar p+\gamma\bar u \in \HOne,
\end{equation}
 which can be interpreted as a Lagrange multiplier for the inequality constraints in
 \eqref{eq:ControlConstraints},
 plays a key role in the error analysis in Section~\ref{sec:Analysis}.
\par
 We can write
\begin{equation}\label{eq:Hahn}
  \lambda=\lambda_1+\lambda_2,
\end{equation}
 where
\begin{equation}\label{eq:lambdaDefinitions}
  \lambda_1=\max(\lambda,0)\geq0 \quad \text{and}\quad
   \lambda_2=\min(\lambda,0)\leq0,
\end{equation}
 and, in view of \eqref{eq:OptimalControlBdd}, \eqref{eq:Multiplier}
  (and $\gamma\leq1$),
\begin{equation}\label{eq:LambdaHOneBounds}
  |\lambda_1|_\HOne,|\lambda_2|_\HOne\leq |\lambda|_\HOne\leq |\bar p|_\HOne
  +\max\big(|\phi_1|_\HOne,|\phi_2|_\HOne,|\bar p|_\HOne\big).
\end{equation}
\par
 From \eqref{eq:ControlCompatibility} and \eqref{eq:baruAndbarp} we have
\begin{equation*}
  \bar u=\begin{cases}
    \phi_2&\qquad \text{if}\; -(\bar p/\gamma)\geq\phi_2\\[4pt]
    -(\bar p/\gamma)&\qquad\text{if}\; \phi_1<-(\bar p/\gamma)<\phi_2\\[4pt]
    \phi_1&\qquad \text{if}\; -(\bar p/\gamma)\leq\phi_1
  \end{cases},
\end{equation*}
 which implies through \eqref{eq:Multiplier} and \eqref{eq:lambdaDefinitions}
  the following complementarity conditions:
\begin{equation}\label{eq:Complementarity}
  \int_\O\lambda_1(\bar u-\phi_1)dx=0=\int_\O \lambda_2(\bar u-\phi_2)dx.
\end{equation}
\begin{remark}\label{rem:Bounds}\rm
 In view of \eqref{eq:baruLTBdd}, \eqref{eq:CSharp}--\eqref{eq:OptimalControlBdd},
  and \eqref{eq:LambdaHOneBounds},
 $\|\bar u\|_\LT$, $|\bar u|_\HOne$, $|\lambda_1|_\HOne$ and $|\lambda_2|_\HOne$
 are bounded by constants
 that only depend on $\|y_d\|_\LT$, $\|f\|_\LT$,
 $\|\phi_1\|_\HOne$, $\|\phi_2\|_\HOne$, $\alpha^{-1}$
  and $\gamma^{-1}$.
\end{remark}
\section{The Approximation Problem}\label{sec:Approximation}
 Let $V_*$ (resp. $W_\dag$) be a closed subspace of $\HOnez$ (resp., $\LT$).
 The approximation problem for \eqref{eq:OCP} is to find
\begin{equation}\label{eq:Discrete}
  (\bysd,\busd)=\argmin_{(y_*,u_\dag)\in \Ksd}
  \frac12\Big[\|y_*-y_d\|_\LT^2+\gamma\|u_\dag\|_\LT^2\Big],
\end{equation}
 where $(y_*,u_\dag)$ belongs to $\Ksd\subset \DES$ if and only if
\begin{equation}\label{eq:DiscretePDEConstraint}
  a(y_*,z_*)=\int_\O (f+u_\dag) z_*dx\qquad\forall\,z_*\in V_*
\end{equation}
 and
\begin{equation}\label{eq:DiscreteControlConstraint}
  Q_\dag \phi_1\leq u_\dag \leq Q_\dag \phi_2\;\;\text{in}\;\;\O.
\end{equation}
Here $Q_\dag:\LT\longrightarrow W_\dag$ is the $L_2$ projection
operator and we assume that
\begin{equation}\label{eq:Order}
   Q_\dag v\geq0 \quad\text{if}\quad v\geq0.
\end{equation}
\par
 Again by the classical theory the minimization problem defined by
 \eqref{eq:Discrete}--\eqref{eq:DiscreteControlConstraint} has a unique
 solution $(\bysd,\busd)\in \Ksd$
 characterized by the first order optimality condition
\begin{equation}\label{eq:DiscreteOptimalityCondition}
  \int_\O (\bysd-y_d)(y_*-\bysd)dx+\gamma\int_O\busd(u_\dag-\busd)dx\geq 0
   \qquad\forall\, (y_*,u_\dag)\in \Ksd.
\end{equation}
\par
 Let $\bpsd\in V_*$ be defined by
\begin{equation}\label{eq:DiscreteAdjointState}
  a(q_*,\bpsd)=\int_\O (\bysd-y_d)q_*dx\qquad\forall\,q_*\in V_*.
\end{equation}
\par
 We will provide estimates for $\|\bar y-\bysd\|_\LT$,
 $\|\bar u-\busd\|_\LT$, $\|\bar p-\bpsd\|_\LT$,
 $|\bar y-\bysd|_a$
 and $|\bar p-\bpsd|_a$ in Section~\ref{sec:Analysis}, where
\begin{equation*}
  |v|_a=\sqrt{a(v,v)} \qquad\forall\,v\in\HOne.
\end{equation*}
\par
 The  simple result below is useful for the analysis of the approximation problem.
 \begin{lemma}\label{lem:RoughCoefficient}
   Let $g\in\LT$ and $v_*\in V_*$ satisfy
 \begin{equation}\label{eq:BVP}
   a(v_*,w_*)=\int_\O gw_* dx\qquad\forall\,w_*\in V_*.
 \end{equation}
  We have
 \begin{align}
    \|v_*\|_\LT&\leq (\CPF^2/\alpha)\|g\|_\LT,\label{eq:BVPLTEst}\\
    |v_*|_a&\leq (\CPF/\sqrt{\alpha})\|g\|_\LT.\label{eq:BVPHOneEst}
 \end{align}
 \end{lemma}
 \begin{proof} The estimate \eqref{eq:BVPLTEst} follows from \eqref{eq:BilinearForm},
 \eqref{eq:PF} and \eqref{eq:BVP}:
 \begin{align*}
  \|v_*\|_\LT^2\leq \CPF^2|v|_\HOne^2&\leq(\CPF^2/\alpha)\,a(v_*,v_*)\\
         &=(\CPF^2/\alpha)\int_\O gv_*dx\leq (\CPF^2/\alpha)\|g\|_\LT\|v_*\|_\LT.
 \end{align*}
 \par
  Similarly we have
\begin{align*}
  |v_*|_a^2=a(v_*,v_*)=\int_\O gv_*dx
  &\leq \|g\|_\LT\|v_*\|_\LT\\
  &\leq \|g\|_\LT\CPF|v_*|_\HOne\leq
  \|g\|_\LT(\CPF/\sqrt{\alpha})|v_*|_a
\end{align*}
 by \eqref{eq:BilinearForm}, \eqref{eq:PF} and \eqref{eq:BVP}, which implies
 \eqref{eq:BVPHOneEst}.
 \end{proof}
\section{Error Estimates}\label{sec:Analysis}
 We will derive error estimates in terms of the $L_2$ projection $Q_\dag:\LT\longrightarrow W_\dag$ and the
 Ritz projection
 $R_*:\HOnez\longrightarrow V_*$  defined by
\begin{equation}\label{eq:Ritz}
  a(R_*\zeta,v_*)=a(\zeta,v_*) \qquad \forall\,v_*\in V_*.
\end{equation}
%
\subsection{Estimate for the $L_2$ Errors}\label{subsec:LTErrors}
\begin{theorem} \label{thm:AbstractLTError}
 There exists a positive constant $C_\flat$ depending only on
 $\alpha^{-1}$ and $\gamma^{-1}$ such that
\begin{align}\label{eq:AbstractLTError}
  &\|\bar y-\bysd\|_\LT+\|\bar u-\busd\|_\LT+\|\bar p-\bpsd\|_\LT\notag\\
  &\hspace{30pt}\leq
  C_\flat\big(\|\bar y-R_*\bar y\|_\LT+\|\bar p-R_*\bar p\|_\LT
  +\|\lambda_1-Q_\dag\lambda_1\|_\LT+\|\lambda_2-Q_\dag\lambda_2\|_\LT\\
  &\hspace{100pt}
   +\|\phi_1-Q_\dag \phi_1\|_\LT+\|\phi_2-Q_\dag\phi_2\|_\LT
   +\|\bar u-Q_\dag\bar u\|_\LT\big).\notag
\end{align}
\end{theorem}
\begin{proof} First we note that \eqref{eq:AdjointState},
\eqref{eq:DiscretePDEConstraint} and
 \eqref{eq:Ritz} imply
\begin{align}\label{eq:DiscreteSAS}
  \int_\O (\bar y-y_d)(y_*-\bysd)dx&=a(y_*-\bysd,\bar p)\\
    &=a(y_*-\bysd,R_*\bar p)=\int_\O (u_\dag-\busd)R_*\bar p\,dx
                                    \notag
\end{align}
for all $(y_*,u_\dag)\in \Ksd$.
\par
 Let $(\tud,\tys)\in \DES$ be defined by
\begin{equation}\label{eq:tud}
  \tud=Q_\dag\bar u
\end{equation}
 and
 \begin{equation}\label{eq:tys}
  a(\tys,z_*)=\int_\O (f+\tud)z_* dx\qquad\forall\,z_*\in V_*.
\end{equation}
 Then $\tud$ satisfies the constraint \eqref{eq:DiscreteControlConstraint}
  by \eqref{eq:ControlConstraints}
 and \eqref{eq:Order}, and hence
 $(\tys,\tud)$ belongs to $\Ksd$.
\par
 We have
\begin{align}\label{eq:LTEst1}
  &\|\bar y-\bysd\|_\LT^2+\gamma\|\bar u-\busd\|_\LT^2\notag\notag\\
  &\hspace{40pt}=\int_\O (\bar y-\bysd)(\bar y-\tys)dx
  +\gamma\int_\O (\bar u-\busd)(\bar u-\tud)dx\\
  &\hspace{70pt}+\int_\O (\bar y-\bysd)(\tys-\bysd)dx
  +\gamma\int_\O(\bar u-\busd)(\tud-\busd)dx.\notag
\end{align}
\par
 Using \eqref{eq:Multiplier}, \eqref{eq:DiscreteOptimalityCondition}
 and \eqref{eq:DiscreteSAS}, we
 find
\begin{align}\label{eq:LTEst2}
  &\int_\O (\bar y-\bysd)(\tys-\bysd)dx+
  \gamma\int_\O(\bar u-\busd)(\tud-\busd)dx\notag\\
  &\hspace{50pt}=\int_\O \bar y(\tys-\bysd)dx+
  \gamma\int_\O \bar u(\tud-\busd)dx\notag\\
  &\hspace{80pt}
   -\int_\O \bysd(\tys-\bysd)dx-\gamma\int_\O \busd(\tud-\busd)dx\\
   &\hspace{50pt}
   \leq \int_\O (\bar y-y_d)(\tys-\bysd)dx+
   \gamma\int_\O \bar u(\tud-\busd)dx\notag\\
   &\hspace{50pt}=\int_\O (R_*\bar p+\gamma\bar u)(\tud-\busd)dx\notag\\
   &\hspace{50pt}=\int_\O \lambda(\tud-\busd)dx+
   \int_\O (R_*\bar p-\bar p)(\tud-\busd)dx,\notag
\end{align}
 and
\begin{equation}\label{eq:LTEst3}
  \int_\O (R_*\bar p-\bar p)(\tud-\busd)dx\leq
   \|\bar p-R_*\bar p\|_\LT\big(\|Q_\dag\bar u-u\|_\LT+\|u-\busd\|_\LT\big)
\end{equation}
 by \eqref{eq:tud}, the Cauchy-Schwarz inequality and the triangle inequality.
\par
 We can estimate the first term on the right-hand side of \eqref{eq:LTEst2} by
 \eqref{eq:Hahn}, \eqref{eq:Complementarity},
 \eqref{eq:DiscreteControlConstraint} and \eqref{eq:tud} as follows.
\begin{align*}
  \int_\O \lambda(\tud-\busd)dx&=\int_\O \lambda_1(\tud-\busd)dx
  +\int_\O\lambda_2(\tud-\busd)dx\notag\\
  &=\int_\O \lambda_1(Q_\dag \bar u-\bar u)dx
  +\int_\O\lambda_2(Q_\dag \bar u-\bar u)dx\notag\\
  &\hspace{30pt}+\int_\O \lambda_1(\bar u-\phi_1)dx
  +\int_\O \lambda_2(\bar u-\phi_2)dx\notag\\
  &\hspace{50pt}
  +\int_\O\lambda_1(\phi_1-Q_\dag\phi_1)dx+
  \int_\O\lambda_2(\phi_2-Q_\dag\phi_2)dx\notag\\
  &\hspace{70pt}
  +\int_\O\lambda_1(Q_\dag\phi_1-\busd)dx+
  \int_\O\lambda_2(Q_\dag\phi_2-\busd)dx\notag\\
  &\leq \int_\O \lambda_1(Q_\dag \bar u-\bar u)dx
  +\int_\O\lambda_2(Q_\dag \bar u-\bar u)dx\notag\\
  &\hspace{30pt}+\int_\O\lambda_1(\phi_1-Q_\dag\phi_1)dx
  +\int_\O\lambda_2(\phi_2-Q_\dag\phi_2)dx\notag\\
  &=\int_\O (\lambda_1-Q_\dag\lambda_1)(Q_\dag \bar u-\bar u)dx
  +\int_\O(\lambda_2-Q_\dag\lambda_2)(Q_\dag \bar u-\bar u)dx\notag\\
  &\hspace{30pt}+\int_\O(\lambda_1-Q_\dag\lambda_1)(\phi_1-Q_\dag\phi_1)dx
  +\int_\O(\lambda_2-Q_\dag\lambda_2)(\phi_2-Q_\dag\phi_2)dx,\notag
\end{align*}
 which  implies
\begin{align}\label{eq:LTEst4}
  \int_\O\lambda(\tud-\busd)dx&\leq \big(\|Q_\dag\bar u-\bar u\|_\LT+
  \|\phi_1-Q_\dag\phi_1\|_\LT+\|\phi_2-Q_\dag\phi_2\|_\LT\big) \\
  &\hspace{40pt}\times \big(\|\lambda_1-Q_\dag\lambda_1\|_\LT+
  \|\lambda_2-Q_\dag\lambda_2\|_\LT\big).\notag
\end{align}
\par
 Putting \eqref{eq:tud} and \eqref{eq:LTEst1}--\eqref{eq:LTEst4}
 together, we arrive at
 the estimate
\begin{align*}
  &\|\bar y-\bysd\|_\LT^2+\gamma\|\bar u-\busd\|_\LT^2\notag\\
  &\hspace{30pt}\leq
   \|\bar y-\bysd\|_\LT\|\bar y-\tys\|_\LT+
   \gamma\|\bar u-\busd\|_\LT\|\bar u-Q_\dag\bar u\|_\LT\notag\\
   &\hspace{60pt}+\|\bar p-R_*\bar p\|_\LT
   \big(\|Q_\dag\bar u-u\|_\LT+\|u-\busd\|_\LT\big)\\
   &\hspace{80pt}+\big(\|Q_\dag\bar u-\bar u\|_\LT+
  \|\phi_1-Q_\dag\phi_1\|_\LT+\|\phi_2-Q_\dag\phi_2\|_\LT\big)\notag \\
  &\hspace{120pt}\times \big(\|\lambda_1-Q_\dag\lambda_1\|_\LT+
  \|\lambda_2-Q_\dag\lambda_2\|_\LT\big),
  \notag
  \end{align*}
 which together with the inequality of arithmetic and geometric means implies
\begin{align}\label{eq:LTEst5}
 &\|\bar y-\bysd\|_\LT^2+\gamma\|\bar u-\busd\|_\LT^2\notag\\
 &\hspace{30pt}
 \leq C\big(\|\bar y-\tys\|_\LT^2+\|\bar u-Q_\dag\bar u\|_\LT^2+\gamma^{-1}\|\bar p-R_*\bar p\|_\LT^2
 +\|\lambda_1-Q_\dag\lambda_1\|_\LT^2\\
 &\hspace{60pt}
   +\|\lambda_2-Q_\dag\lambda_2\|_\LT^2+
    \|\phi_1-Q_\dag\phi_1\|_\LT^2+\|\phi_2-Q_\dag\phi_2\|_\LT^2\big),
   \notag
\end{align}
 where $C$ is a universal positive constant.
\par
 Note that \eqref{eq:PDEConstraint}, {\eqref{eq:Ritz}}, \eqref{eq:tud} and \eqref{eq:tys} imply
\begin{equation*}
 a(R_*\bar y-\tys,z_*) = a(\bar y-\tys,z_*)=\int_\O (\bar u-\tud)z_*dx
   =\int_\O (\bar u-Q_\dag\bar u)z_*dx\qquad\forall\,z_*\in V_*
\end{equation*}
 and hence
\begin{equation*}
  \|{ R_*\bar y}-\tys\|_\LT\leq (\CPF^2/\alpha)\|\bar u-Q_\dag\bar u\|_\LT
\end{equation*}
 by  Lemma~\ref{lem:RoughCoefficient}.  { Therefore we have
\begin{equation}\label{eq:LTEst6}
  \|\bar y-\tys\|_\LT\leq \|\bar y-R_*\bar y\|_\LT+(\CPF^2/\alpha)\|\bar u-Q_\dag\bar u\|_\LT.
\end{equation}
}
\par
 Similarly \eqref{eq:AdjointState}, \eqref{eq:DiscreteAdjointState}
 and \eqref{eq:Ritz} imply
\begin{equation}\label{eq:AdjointStateRelation}
  a(q_*,R_*\bar p-\bpsd)=a(q_*,\bar p-\bpsd)=
  \int_\O (\bar y-\bysd)q_*dx\qquad\forall\,q_*\in V_*,
\end{equation}
 and hence
\begin{equation}\label{eq:Ritzbarp}
  \|R_*\bar p-\bpsd\|_\LT\leq (\CPF^2/\alpha)\|\bar y-\bysd\|_\LT
\end{equation}
 by Lemma~\ref{lem:RoughCoefficient}.  Consequently we have
\begin{equation}\label{eq:LTEst7}
  \|\bar p-\bpsd\|_\LT
  \leq
   \|\bar p-R_*\bar p\|_\LT+(\CPF^2/\alpha)\|\bar y-\bysd\|_\LT
\end{equation}
 by \eqref{eq:Ritzbarp} and the triangle inequality.
\par
 The estimate \eqref{eq:AbstractLTError} follows from \eqref{eq:LTEst5},
 \eqref{eq:LTEst6} and \eqref{eq:LTEst7}.
\end{proof}
\par
 The following result shows that \eqref{eq:AbstractLTError} is
 sharp up to the terms involving $Q_\dag$.
\begin{theorem}\label{thm:Tight}
 There exists a positive constant $C_\natural$ depending only on
 $\alpha^{-1}$ such that
\begin{equation}\label{eq:Tight}
  \|\bar y-R_*\bar y\|_\LT+\|\bar p-R_*\bar p\|_\LT
  \leq C_\natural\big(\|\bar y-\bysd\|_\LT+\|\bar u-\busd\|_\LT
  +\|\bar p-\bpsd\|_\LT\big).
\end{equation}
\begin{proof} We have
\begin{equation}\label{eq:Tight1}
  \|\bar y-R_*\bar y\|_\LT\leq \|\bar y-\bysd\|_\LT+\|\bysd-R_*\bar y\|_\LT,
\end{equation}
 and, in view of \eqref{eq:PDEConstraint},
 \eqref{eq:DiscretePDEConstraint} and \eqref{eq:Ritz},
\begin{equation}\label{eq:StateRelation}
  a(R_*\bar y-\bysd,z_*)=a(\bar y-\bysd,z_*)=
  \int_\O (\bar u-\busd)z_*dx\qquad\forall\,z_*\in V_*,
\end{equation}
 which implies through Lemma~\ref{lem:RoughCoefficient}
\begin{equation}\label{eq:Tight2}
  \|\bysd-R_*\bar y\|_\LT\leq (\CPF^2/\alpha)\|\bar u-\busd\|_\LT.
\end{equation}
\par
 Similarly we have
\begin{equation}\label{eq:Tight3}
  \|\bar p-R_*\bar p\|_\LT\leq \|\bar p-\bpsd\|_\LT
  +(\CPF^2/\alpha)\|\bar y-\bysd\|_\LT
\end{equation}
 by \eqref{eq:Ritzbarp} and the triangle inequality.
 \par
  The estimate \eqref{eq:Tight} follows from \eqref{eq:Tight1},
  \eqref{eq:Tight2} and \eqref{eq:Tight3}.
\end{proof}
\end{theorem}
\subsection{Estimate for the Energy Errors}\label{subsec:EnergyErrors}
\begin{theorem} \label{thm:AbstractEnergyError} There exists a positive
constant $C_\maltese$ depending only on $\alpha^{-1}$ and $\gamma^{-1}$ such that
\begin{align}\label{eq:AbstractEnergyError}
   &|\bar y-\bysd|_a+|\bar p-\bpsd|_a\notag\\
   &\hspace{40pt}\leq C_\maltese \big(|\bar y-R_*\bar y|_a+|\bar p-R_*\bar p|_a
   +\|\lambda_1-Q_\dag\lambda_1\|_\LT+\|\lambda_2-Q_\dag\lambda_2\|_\LT\\
   &\hspace{120pt}+\|\phi_1-Q_\dag\phi_1\|_\LT+\|\phi_2-Q_\dag\phi_2\|_\LT
   +\|\bar u-Q_\dag\bar u\|_\LT\big).
   \notag
\end{align}
\end{theorem}
\begin{proof}
  We begin with a triangle inequality
\begin{equation}\label{eq:Energy1}
  |\bar y-\bysd|_a+|\bar p-\bpsd|_a\leq |\bar y-R_*\bar y|_a+|R_*\bar y-\bysd|_a+
  |\bar p-R_*\bar p|_a+|R_*\bar p-\bpsd|_a.
\end{equation}
\par
 From \eqref{eq:StateRelation} we obtain
\begin{equation}\label{eq:Energy2}
  |R_*\bar y-\bysd|_a\leq (\CPF/\sqrt{\alpha})\|\bar u-\busd\|_\LT
\end{equation}
 by Lemma~\ref{lem:RoughCoefficient}.
\par
 Similarly we have
\begin{equation}\label{eq:Energy3}
  |R_*\bar p-\bpsd|_a\leq (\CPF/\sqrt{\alpha})\|\bar y-\bysd\|_\LT
\end{equation}
 by \eqref{eq:AdjointStateRelation} and Lemma~\ref{lem:RoughCoefficient}.
\par
 Finally the Poincar\'e-Friedrichs inequality \eqref{eq:PF} and \eqref{eq:BilinearForm} imply
\begin{equation}\label{eq:Energy4}
  \|\bar y-R_*\bar y\|_\LT+\|\bar p-R_*\bar p\|_\LT\leq
  (\CPF/\sqrt\alpha)\big(|\bar y-R_*\bar y|_a+|\bar p-R_*\bar p|_a\big),
\end{equation}
 and the estimate \eqref{eq:AbstractEnergyError} follows from
 Theorem~\ref{thm:AbstractLTError}
 and
 \eqref{eq:Energy1}--\eqref{eq:Energy4}.
\end{proof}
\section{Applications}\label{sec:Applications}
 We can apply the error estimates in Section~\ref{sec:Analysis}
 to standard finite element methods and multiscale finite element methods.
\subsection{Standard Finite Element Methods}\label{subsec:Classical}
 We assume that the bilinear form $a(\cdot,\cdot)$ is given by
\begin{equation}\label{eq:ClasicalBilinearForm}
  a(y,z)= \int_\O \big[A(x)\nabla y\cdot\nabla z + c(x)yz\big]dx,
\end{equation}
 where the nonnegative function $c(x)$  and the $d\times d$
 symmetric matrix function $A(x)$
   are sufficiently smooth, and there exists a positive
    constant $\mu$ such that
\begin{equation*}
  \xi^t A(x)\xi\geq \mu|\xi|^2 \qquad\forall\,x\in\O,\;\xi\in\R^d.
  \end{equation*}
\par
 We can take $V_*=V_h\subset \HOnez$ to be the $P_1$ Lagrange finite element space
 (cf. \cite{Ciarlet:1978:FEM,BScott:2008:FEM}) associated with
 a regular triangulation $\cT_h$ of $\O$, and $W_\dag=W_\rho$ to be the space
 of piecewise constant functions
 associated with a regular triangulation $\cT_\rho$ of $\O$.
 The optimal state (resp., optimal control
 and adjoint state) is denoted by $\bar y_{h,\rho}$ (resp.,
 $\bar u_{h,\rho}$ and $\bar p_{h,\rho}$).
\par
 For simplicity, we assume $\O$ is convex.
 It is known that $\bar y$ and $\bar p$ belong to $H^2(\O)$
 (cf. \cite{Grisvard:1985:EPN,Dauge:1988:EBV,MR:2010:Polyhedral})
  and we have the following estimates
 (cf. \cite{Ciarlet:1978:FEM,BScott:2008:FEM}) for the
 Ritz projection operator
 $R_h:\HOnez\longrightarrow V_h$:
\begin{alignat}{3}
  |\zeta-R_h\zeta|_a&
  \leq C_1h|\zeta|_{H^2(\O)}&\qquad&\forall\,
  \zeta\in H^2(\O)\cap \HOnez,\label{eq:StandardRitzEnergy}\\
  \|\zeta-R_h\zeta\|_\LT&\leq
  C_1h^{2}|\zeta|_{H^2(\O)}&\qquad&\forall\,\zeta\in H^2(\O)\cap \HOnez,
  \label{eq:StandardRitzLT}
\end{alignat}
 where the positive constant $C_1$ depends only on the
 coefficients in \eqref{eq:ClasicalBilinearForm}
 and the shape regularity of $\cT_h$.
\par
 The  $L_2$ projection $Q_\rho:\LT\longrightarrow W_\rho$
 satisfies \eqref{eq:Order} and we have a standard
 error estimate
 (cf. \cite{Ciarlet:1978:FEM,BScott:2008:FEM}):
\begin{equation}\label{eq:StandardLTProjection}
  \|\zeta-Q_\rho\zeta\|_\LT\leq C_2h|\zeta|_\HOne\qquad\forall\,\zeta\in \HOne,
\end{equation}
 where the positive constant $C_2$ depends only on
 the shape regularity of $\cT_\rho$.
\par
 It follows from Remark~\ref{rem:Bounds}, Theorem~\ref{thm:AbstractLTError},
 Theorem~\ref{thm:AbstractEnergyError},
 and \eqref{eq:StandardRitzEnergy}--\eqref{eq:StandardLTProjection} that
\begin{align}
  \|\bar y-\bar y_{h,\rho}\|_\LT+\|\bar u-\bar u_{h,\rho}\|_\LT+
   \|\bar p-\bar p_{h,\rho}\|_\LT&\leq C(h^2+\rho),\label{eq:Falk1}\\
   |\bar y-\bar y_{h,\rho}|_a+|\bar p-\bar p_{h,\rho}|_a&\leq C(h+\rho),
   \label{eq:Falk2}
\end{align}
 where the positive constant $C$ is independent of $h$ and $\rho$, and we have
 recovered the error estimates in \cite{Falk:1973:Control}
 for a convex $\O$.
\par
 We can also take $W_\dag$ to be $\LT$, which is the variational
 discretization concept in
 \cite{Hinze:2005:Control}.  In this case $Q_\dag$ is the identity
 map on $\LT$ so that \eqref{eq:Order}
  is satisfied trivially and we denote the
 optimal state (resp., optimal control and adjoint state) by $\bar y_h$
  (resp., $\bar u_h$ and $\bar p_h$).  The
 estimates \eqref{eq:Falk1} and \eqref{eq:Falk2} become
\begin{align}
  \|\bar y-\bar y_h\|_\LT+\|\bar u-\bar u_h\|_\LT+\|\bar p-\bar p_h\|_\LT&\leq Ch^2,
  \label{eq:Hinze1}\\
  |\bar y-\bar y_h|_a+|\bar p-\bar p_h|_a&\leq Ch,\label{eq:Hinze2}
\end{align}
 where $C$ is independent of $h$, and we have recovered the result in
 \cite{Hinze:2005:Control}.
\begin{remark}\label{rem:GradedMesh}\rm
  The estimates \eqref{eq:Falk1}--\eqref{eq:Hinze2} also hold for
  a general $\O$ provided
  the triangulations $\cT_h$ and $\cT_\rho$ are properly graded around the
  singular parts of $\p\O$
  (cf. \cite{Li:2022:Book}).
\end{remark}
\subsection{Multiscale Finite Element Methods}\label{subsec:Multiscale}
 Under assumption \eqref{eq:BilinearForm}, the optimal state $\bar y$ and
 adjoint state $\bar p$ belong
 to $\HOnez$ and we cannot claim any additional regularity.
 \par  If we take
 $V_*=V_h\subset \HOnez$ to be
 the $P_1$ finite element space associated with $\O$ and $W_\dag=W_\rho$ to be
 the space of piecewise constant functions
 associated with $\cT_\rho$, then Theorem~\ref{thm:AbstractLTError} implies
\begin{align*}
  &\lim_{h,\rho\downarrow0}\big(\|\bar y-\bar y_{h,\rho}\|_\LT+
  \|\bar u-\bar u_{h,\rho}\|_\LT+
   \|\bar p-\bar p_{h,\rho}\|_\LT\big)\\
   &\hspace{40pt}\leq \lim_{h\downarrow 0}\big(\|\bar y-R_h\bar y\|_\LT+
   \|\bar p-R_h\bar p\|_\LT\big)\\
   &\hspace{70pt}+\lim_{\rho\downarrow0}\big(\|\lambda_1-Q_\rho\lambda_1\|_\LT
   +\|\lambda_1-Q_\rho\lambda_1\|_\LT+\|\phi_1-Q_\rho\phi_1\|_\LT\\
  &\hspace{120pt}\|\phi_2-Q_\rho\phi_2\|_\LT+\|\bar u-Q_\rho \bar u\|_\LT\big)\\
   &\hspace{40pt}=0.
\end{align*}
 Therefore this standard finite element method converges, but the convergence in
 $h$ can be arbitrarily
 slow (cf. \cite{BO:2000:Bad}), and an accurate approximation of
 $(\bar y,\bar u,\bar p)$ will require
 a very small mesh size $h$.
\par
 We can remedy this slow convergence by taking $V_*$ to be a multiscale
  finite element space.  For example  we can take $V_*$ to be the rough
  polyharmonic space
  $V_H^{rps}$ in \cite{OZB:2014:Homogenization,LZZ:2021:Polyharmonic}
  associated with a triangulation
  $\cT_H$ and $W_\dag=W_\rho$ remains the space of piecewise constant
  functions associated with a triangulation $\cT_\rho$.  The optimal
  state (resp., optimal
   control and adjoint state) is denoted by $\bar y_{H,\rho}^{rps}$
    (resp., $\bar u_{H,\rho}^{rps}$ and
   $\bar p_{H,\rho}^{rps}$).
 \par
  Let $\zeta\in\HOnez$ satisfy
\begin{equation}\label{eq:zetaEq}
  a(\zeta,v_{\sss H})=\int_\O gv_{\sss H}dx \qquad\forall\,v_{\sss H}\in V_H^{rps},
\end{equation}
 where $g\in\LT$.
 Then we have, by \eqref{eq:BilinearForm}, \eqref{eq:PF} and the
 estimates in \cite{OZB:2014:Homogenization,LZZ:2021:Polyharmonic} ,
\begin{equation}\label{eq:RPSEstimate}
 \|\zeta-R_H^{rps}\zeta\|_\LT\leq (\CPF/\sqrt{\alpha})
  |\zeta-R_H^{rps}\zeta|_a \leq C_3 H \|g\|_\LT,
\end{equation}
 where $R_H^{rps}:\HOnez\longrightarrow V_H^{rps}$ is the Ritz projection operator
 and the positive constant $C_3$ depends only on the shape regularity
 of $\cT_H$ and $\alpha^{-1}$.
\par
 It follows from Remark~\ref{rem:Bounds},
 Theorem~\ref{thm:AbstractLTError}, Theorem~\ref{thm:AbstractEnergyError}
 and \eqref{eq:RPSEstimate} that
\begin{align}\label{eq:MSEstimate}
  &\|\bar y-\bar y_{H,\rho}^{rps}\|_\LT+\|\bar u-\bar u_{H,\rho}^{rps}\|_\LT
  +\|\bar p-\bar p_{H,\rho}^{rps}\|_\LT+
  |\bar y-\bar y_{H,\rho}^{rps}|_a+|\bar p-\bar p_{H,\rho}^{rps}|_a\\
  &\hspace{40pt}\leq C_\diamond (H+\rho),\notag
\end{align}
 where the positive constant $C_\diamond$ depends only on
 $\|y_d\|_\LT$, $\|f\|_\LT$,
 $\|\phi_1\|_\HOne$, $\|\phi_2\|_\HOne$, $\alpha^{-1}$,
  $\gamma^{-1}$, and the shape regularities of $\cT_H$ and $\cT_\rho$,
   and we have recovered the
  results in \cite{CLZZ:2023:MultiscaleCC}.
\par
 We can also take $V_*$ to be the constraint energy minimizing generalized multiscale
 finite element space $V_H^{gms}$ in \cite{CEL:2018:Multiscale}
 associated with a triangulation $\cT_H$.  In this case the function
 $\zeta\in\HOnez$ defined by
 \eqref{eq:zetaEq} satisfies
\begin{equation*}
  \|\zeta-R_H^{gms}\zeta\|_\LT\leq (\CPF/\sqrt{\alpha})
  |\zeta-R_H^{gms}\zeta|_a\leq C_4H\|g\|_\LT,
\end{equation*}
 where the positive constant $C_4$ depends only on $\alpha^{-1}$,
  the shape regularity of $\cT_H$ and
 $\Lambda^{-1}$.  ($\Lambda$ is a spectral parameter used in the
 construction of the multiscale
 finite element space $V_H^{gms}$.)
 Therefore \eqref{eq:MSEstimate} also holds for the approximate solution
 $(\bar y_{H,\rho}^{gms},\bar u_{H,\rho}^{gms},\bar p_{H,\rho}^{gms})$
 obtained by this multiscale
 finite element method where $C_\diamond$ depends also on $\Lambda^{-1}$.
  This is the result in \cite{AC:2022:DD26} in the case where $\cT_\rho=\cT_h$.
\par
 Finally we can take $V_*$ to be the local orthogonal decomposition multiscale
  finite element spaces $V_H^{lod}$
  in \cite{MP:2014:LOD,MP:2021:LOD,BGS:2021:LOD} associated with
 a triangulation $\cT_H$ that incorporates information from
 a standard finite element space
 $V_h$ associated with a refinement $\cT_h$ of $\cT_H$.
 We denote the optimal state (resp., optimal
   control and adjoint state) by $\bar y_{H,\rho}^{lod}$
   (resp., $\bar u_{H,\rho}^{lod}$ and
   $\bar p_{H,\rho}^{lod}$), and the Ritz projection operator from
   $\HOnez$ to $V_H^{lod}$ is denoted by
   $R_H^{lod}$.
\par
 Let $v_h\in V_h$ and $v_{\sss H}\in V_H^{lod}$ satisfy
\begin{alignat*}{3}
  a(v_h,w_h)&=\int_\O gw_h dx&\qquad&\forall\,w_h\in V_h,\\
  a(v_{\sss H},w_{\sss H})&=\int_\O gw_{\sss H}dx
  &\qquad&\forall\,w_{\sss H}\in V_H^{lod}.
\end{alignat*}
 Then we have, by the results in \cite{MP:2014:LOD} and \cite{BGS:2021:LOD},
\begin{equation}\label{eq:LODEstimates}
  |v_h-v_{\sss H}|_a\leq C_5 H\|g\|_\LT \quad\text{and}\quad
  \|v_h-v_{\sss H}\|_\LT\leq C_5H^2\|g\|_\LT,
\end{equation}
 where the positive constant $C_5$ depends only on
  $\alpha^{-1}$ and the shape regularity
 of $\cT_H$.
\par
 According to Remark~\ref{rem:Bounds}, Theorem~\ref{thm:AbstractLTError} and
 \eqref{eq:StandardLTProjection}, we have
\begin{align}\label{eq:LODEst1}
  &\|\bar y- \bar y_{H,\rho}^{lod}\|_\LT+
  \|\bar u-\bar u_{H,\rho}^{lod}\|_\LT+\|\bar p-\bar p_{H,\rho}^{lod}\|_\LT\\
  &\hspace{60pt}\leq
  C_6 \big(\|\bar y-R_H^{lod}\bar y\|_\LT+
  \|\bar p-R_H^{lod}\bar p\|_\LT+\rho\big),\notag
\end{align}
 where the positive constant $C_6$ only depends on
 $\|y_d\|_\LT$, $\|f\|_\LT$,
 $\|\phi_1\|_\HOne$, $\|\phi_2\|_\HOne$, $\alpha^{-1}$,
  $\gamma^{-1}$, and the shape regularity of $\cT_\rho$.
 \par
  Let $(\bar v_{h,\rho},\bar u_{h,\rho},\bar p_{h,\rho})$ be the solution of
  \eqref{eq:Discrete} based on $V_*=V_h$ and $W_\dag=W_\rho$.  Then we have
\begin{alignat*}{3}
  a(R_h\bar y,z_h)&=\int_\O (f+\bar u) z_h dx&\qquad&\forall\,z_h\in V_h,\\
  a(R_H^{lod}\bar y,z_{\sss H})&=\int_\O (f+\bar u)z_{\sss H}dx&\qquad&\forall\,z_{\sss H}\in V_H^{lod},
\end{alignat*}
 by \eqref{eq:PDEConstraint} and \eqref{eq:Ritz}, which together with \eqref{eq:LODEstimates}
 imply
\begin{equation}\label{eq:LODEst2}
  \|R_h\bar y-R_H^{lod}\bar y\|_\LT\leq C_5H^2\|f+\bar u\|_\LT.
\end{equation}
\par
 Similarly we have
\begin{equation}\label{eq:LODEst3}
  \|R_h\bar p-R_H^{lod}\bar p\|_\LT\leq C_5H^2\|\bar y-y_d\|_\LT
\end{equation}
 by \eqref{eq:AdjointState}, \eqref{eq:Ritz} and \eqref{eq:LODEstimates}.
\par
 Putting Theorem~\ref{thm:Tight} and \eqref{eq:LODEst1}--\eqref{eq:LODEst3}
  together we arrive at
\begin{align}\label{eq:LODLT}
   &\|\bar y- \bar y_{H,\rho}^{lod}\|_\LT+\|\bar u-\bar u_{H,\rho}^{lod}\|_\LT
   +\|\bar p-\bar p_{H,\rho}^{lod}\|_\LT
   \notag\\
  &\hspace{60pt}\leq
  C_6 \big(\|\bar y-R_h\bar y\|_\LT+\|R_h\bar y-R_H^{lod}\bar y\|_\LT\notag\\
  &\hspace{100pt}+
   \|\bar p-R_h\bar y\|_\LT+\|R_h\bar y-R_H^{lod}\bar p\|_\LT+\rho\big)\notag\\
   &\hspace{60pt}\leq C_6\big(\|\bar y-R_h\bar y\|_\LT+\|\bar p-R_h\bar p\|_\LT
   \\
   &\hspace{100pt}+ C_5H^2\|f+\bar u\|_\LT+C_5H^2\|\bar y-y_d\|_\LT+\rho\big)\notag\\
   &\hspace{60pt}\leq
    C_6C_\sharp\big(\|\bar y-\bar y_{h,\rho}\|_\LT+\|\bar u-\bar u_{h,\rho}\|_\LT+
     \|\bar p-\bar p_{h,\rho}\|_\LT\big)\notag\\
     &\hspace{100pt}+ C_6C_5H^2\big(\|f+\bar u\|_\LT+\|\bar y-y_d\|_\LT\big)+C_6\rho.\notag
\end{align}
\par
 Similarly, we have by Theorem~\ref{thm:AbstractEnergyError}
\begin{align}\label{eq:LODEnergy}
  |\bar y-\bar y_{H,\rho}^{lod}|_a+|\bar p-\bar p_{H,\rho}^{lod}|_a
  \leq C_7\big(|\bar y-\bar y_{h,\rho}|_a+|\bar p-\bar p_{h,\rho}|_a+H+\rho\big),
\end{align}
 where the positive constant $C_7$ only depends on
 $\|y_d\|_\LT$, $\|f\|_\LT$,
 $\|\phi_1\|_\HOne$, $\|\phi_2\|_\HOne$, $\alpha^{-1}$,
  $\gamma^{-1}$, and the shape regularities of $\cT_H$ and $\cT_\rho$.
\par
 Comparing \eqref{eq:Falk1}--\eqref{eq:Falk2} and \eqref{eq:LODLT}--\eqref{eq:LODEnergy},
 we conclude that up to the error of a fine scale approximation, the
 performance of the local orthogonal decomposition
  multiscale finite element method with rough coefficients and
 on a general $\O$ is identical to the performance of  standard finite element methods for
 a problem with smooth coefficients and on a convex domain.
  Moreover all the constants in the estimates are independent of the
  mesh sizes and the contrast $\beta/\alpha$.
\par
 Numerical results for the local orthogonal decomposition method for \eqref{eq:OCP}
  can be found in \cite{BGS:2024:LOD_CC}.

\section{Conclusions}\label{sec:Conclusions}
 We have developed new abstract error estimates for a model linear-quadratic elliptic
 distributed optimal control problem that reduce the error analysis
  to the properties of the Ritz projection operator for the finite element space
  for the state
   and the $L_2$ projection operator
  for the finite element space for the control.
  They can be applied to standard finite element methods
 for a classical partial differential equation constraint and  multiscale finite element
 methods when the coefficients in the partial differential equation constraint are rough.
 Besides the multiscale finite element methods mentioned in Section~\ref{sec:Applications},
 they can also be applied to many others, such as the ones investigated in
 \cite{HFMQ:1998:VMS,HS:2007:VM,BL:2011:Multiscale,EGH:2013:GMSFEM,MS:2022:Error,MSD:2022:GFEM}.
\par
 For simplicity we have assumed that $a(\cdot,\cdot)$ is symmetric.
  But the estimates in Section~\ref{sec:Analysis}
 can be extended to a
 nonsymmetric $a(\cdot,\cdot)$ by replacing the term $R_*\bar p$ with the term
 $S_*\bar p$, where $S_*:\HOnez\longrightarrow V_*$ is defined by
\begin{equation*}
  a(q_*,S_*\zeta)=a(q_*,\zeta) \qquad\forall\,q_*\in V_*.
\end{equation*}
\par
 Finally we note that error estimates for boundary control problems with
 rough coefficients are still absent.

\end{document}